\documentclass[12pt]{article}
\font\smallit=cmti10 
\font\smalltt=cmtt10 
\font\smallrm=cmr9 

\usepackage{amsmath,latexsym} 

\newcommand{\Z}{\ensuremath{\mathbf Z}}

\newtheorem{theorem}{Theorem} 
\newtheorem{lemma}{Lemma} 
\newtheorem{corollary}{Corollary} 
\newcommand{\bt}{\begin{theorem}} 
\newcommand{\et}{\end{theorem}} 
\newcommand{\bl}{\begin{lemma}} 
\newcommand{\el}{\end{lemma}} 
\newcommand{\bc}{\begin{corollary}} 
\newcommand{\ec}{\end{corollary}} 
\newcommand{\benum}{\begin{enumerate}} 
\newcommand{\eenum}{\end{enumerate}} 
\newcommand{\card}{\ensuremath{\text{card}}} 

\makeatletter 

\renewcommand\section{\@startsection {section}{1}{\z@}%
{-30pt \@plus -1ex \@minus -.2ex}%
{2.3ex \@plus.2ex}%
{\normalfont\normalsize\bfseries}} 

\renewcommand\subsection{\@startsection{subsection}{2}{\z@}%
{-3.25ex\@plus -1ex \@minus -.2ex}%
{1.5ex \@plus .2ex}%
{\normalfont\normalsize\bfseries}} 

\renewcommand{\@seccntformat}[1]{\csname the#1\endcsname. } 

\makeatother

\begin{document} 
\vspace*{-40pt} 
\centerline{\smalltt INTEGERS: \smallrm 
ELECTRONIC JOURNAL OF COMBINATORIAL NUMBER THEORY \smalltt 7 
(2007), \#A01} 
\vskip 40pt 





\begin{center} 
\uppercase{\bf Affine invariants, 
relatively prime sets, and a phi function for subsets of 
$\mathbf{\{1,2,\ldots,n\}}$} 
\vskip 20pt 
{\bf Melvyn B. Nathanson\footnote{This work was supported 
in part by grants from the NSA 
Mathematical Sciences Program and the PSC-CUNY Research Award Program.}}\\ 
{\smallit Lehman College (CUNY), Bronx, New York 10468}\\ 
{\tt melvyn.nathanson@lehman.cuny.edu}\\ 
\end{center} 
\vskip 30pt 
\centerline{\smallit Received: 8/15/06, 
Accepted: 12/29/06, Published: 1/3/07} 
\vskip 30pt 

\centerline{\bf Abstract} 
\noindent 
A nonempty subset $A$ of $\{1,2,\ldots,n\}$ is 
{\em relatively prime} if $\gcd(A)=1.$ Let $f(n)$ and $f_k(n)$ denote, 
respectively, the number of relatively prime subsets and the number of 
relatively prime subsets of cardinality $k$ of $\{1,2,\ldots,n\}$. Let 
$\Phi(n)$ and $\Phi_k(n)$ denote, respectively, the number of nonempty 
subsets and the number of subsets of cardinality $k$ of 
$\{1,2,\ldots,n\}$ such that $\gcd(A)$ is relatively prime to $n$. 
Exact formulas and asymptotic estimates are obtained for these functions. 

\footnotesize 
\noindent 
Subject class: Primary 11A25, 11B05, 11B13, 11B75. 


\normalsize 

\pagestyle{myheadings} 
\markright{\smalltt INTEGERS: \smallrm ELECTRONIC JOURNAL 
OF COMBINATORIAL NUMBER THEORY \smalltt 7 (2007), \#A01\hfill} 

\thispagestyle{empty} 
\baselineskip=15pt 
\vskip 30pt 

\section{Affine Invariants} 
Let $A$ be a set of integers, and let $x$ and $y$ be rational numbers. We define the {\em dilation} $x\ast A = \{xa: a \in A\}$ and the {\em translation} $A + y = \{a+y : a \in A\}$. Sets of integers $A$ and $B$ are {\em affinely equivalent} if 
there exist rational numbers $x \neq 0$ and $y$ such that $B = x\ast A + y$. For example, the sets $A = \{2,8,11,20\}$ and $B = \{-4,10,17,38 \}$ are affinely equivalent, since $B = (7/3)\ast A - 26/3 $, and $A$ and $B$ are both affinely equivalent to the sets $C = \{0,2,3,6\}$ and $D = \{0,3,4,6\}$. Every set with one element is affinely equivalent to $\{0\}$. Every finite set $A$ of integers with more than one element is affinely equivalent to unique sets $C$ and $D$ of nonnegative integers such that $\min(C)=\min(D)=0$, $\gcd(C) = \gcd(D) = 1$, and $D = (-1)\ast C +\max(C).$ 

A function $f(A)$ whose domain is the set $\mathcal{F}(\Z)$ of nonempty finite sets of integers is called an {\em affine invariant} of $\mathcal{F}(\Z)$ if $f(A) = f(B)$ for all affinely equivalent sets $A$ and $B$. For example, if $A+A = \{ a+a' : a,a'\in A\}$ is the sumset of a finite set $A$ of integers, and if $A-A = \{ a-a' : a,a'\in A\}$ is the difference set of the finite set $A$, then $s(A) = \card(A+A)$ and $d(A) = \card(A-A)$ are affine invariants. More generally, let $u_0,u_1,\ldots, u_n$ be integers and $F(x_1,\ldots,x_n) = u_1x_1+\cdots + u_nx_n+u_0$. Define $F(A) = \{u_1a_1+\cdots + u_na_n+u_0 : a_1,\ldots,a_n \in A \text{ for } i=1,\ldots,n\}.$ Then $f(A) = \card(F(A))$ is an affine invariant. 

Let $f(A)$ be a function with domain $\mathcal{F}(\Z)$. A frequent problem in combinatorial number theory is to determine the distribution of values of the function $f(A)$ for sets $A$ in the interval of integers $\{0,1,\ldots,n\}$. 
For example, if $A \subseteq \{0,1,2,\ldots,n\}$, then $1 \leq \card(A+A) \leq 2n+1$. For $\ell = 1,\ldots, 2n+1$, we can ask for the number of nonempty sets $A \subseteq \{0,1,2,\ldots,n\}$ such that $\card(A+A) = \ell$. 
Similarly, if $\emptyset \neq A \subseteq \{0,1,2,\ldots,n\}$ and $\card(A) = k$, then $2k-1 \leq \card(A+A) \leq k(k+1)/2$, and, for $\ell = 2k-1,\ldots,k(k+1)/2$, we can ask for the number of such sets $A$ with $\card(A+A) = \ell$. 
In both cases, there is a redundancy in considering sets that are affinely equivalent, and we might want to count only sets that are pairwise affinely inequivalent.

\section{Relatively Prime Sets} 
A nonempty subset $A$ of $\{1,2,\ldots,n\}$ will be called {\it relatively prime} if the elements of $A$ are relatively prime, that is, if $\gcd(A)=1$. 
Let $f(n)$ denote the number of relatively prime subsets of $\{1,2,\ldots,n\}$. The first 10 values of $f(n)$ are 1, 2, 5, 11, 26, 53, 116, 236, 488, and 983. 
(This is sequence A085945 in Sloane's {\it On-Line Encyclopedia of Integer Sequences}.) 
Let $f_k(n)$ denote the number of relatively prime subsets of $\{1,2,\ldots,n\}$ of cardinality $k$. 
We present exact formulas and asymptotic estimates for $f(n)$ and $f_k(n)$. These estimates imply that almost all finite sets of integers are relatively prime. 

No set of even integers is relatively prime. 
Since there are $2^{[n/2]} - 1$ nonempty subsets of 
$\{2,4,6,\ldots,2[n/2]\}$ and $2^n-1$ nonempty subsets of 
$\{1,2,\ldots,n\}$, we have the upper bound 
\begin{equation} \label{prim:upper} 
f(n) \leq 2^n - 2^{[n/2]}. 
\end{equation} 
Similarly, 
\begin{equation} \label{prim:upper-fkn} 
f_k(n) \leq {n \choose k} - {[n/2] \choose k}. 
\end{equation} 

If $1 \in A$, then $A$ is relatively prime. Since there are $2^{n-1}$ sets $A \subseteq \{1,2,\ldots,n\}$ with $1 \in A$, we have 
\[ 
f(n) \geq 2^{n-1}. 
\] 
Let $n \geq 3.$ If $1 \notin A$ but $2 \in A$ and $3\in A$, then $A$ is relatively prime and so 
\[ 
f(n) \geq 2^{n-1} + 2^{n-3}. 
\] 
Let $n \geq 5.$ If $1 \notin A$ and $3 \notin A$, but $2\in A$ and $5\in A$, then $A$ is relatively prime. If $1 \notin A$ and $2 \notin A$, but $3\in A$ and $5\in A$, then $A$ is relatively prime. Therefore, 
\[ 
f(n) \geq 2^{n-1} + 2^{n-3} + 2\cdot 2^{n-4} = 2^{n-1} + 2^{n-2}. 
\] 
Similarly, 
\[ 
f_k(n) \geq {n-1 \choose k-1} + {n-3 \choose k-2} + 2{n-4 \choose k-2}. 
\]

\section{Exact Formulas and Asymptotic Estimates} 
Let $[x]$ denote the greatest integer less than or equal to $x$. 
If $x \geq 1$ and $n = [x]$, then 
\[ 
\left[ \frac{x}{d} \right] = \left[ \frac{[x]}{d} \right] = \left[ \frac{n}{d} \right] 
\] 
for all positive integers $d$.

Let $F(x)$ be a function defined for $x \geq 1$, and define the function 
\[ 
G(x) = \sum_{1 \leq d \leq x} F\left(\frac{x}{d}\right). 
\] 
In the proof of Theorem~\ref{prim:theorem: f} we use the following version of the M{\" o}bius inversion formula (Nathanson~\cite[Exercise 5 on p. 222]{nath00A}): 
\[ 
F(x) = \sum_{1 \leq d \leq x} \mu(d) G\left(\frac{x}{d}\right). 
\] 

\bt \label{prim:theorem: f} 
For all positive integers $n$, 
\begin{equation} \label{prim:recursion} 
\sum_{d=1}^n f\left(\left[\frac{n}{d}\right]\right) = 2^n-1 
\end{equation} 
and 
\begin{equation} \label{prim:explicit} 
f(n) = \sum_{d=1}^{n} \mu(d) \left(2^{[n/d]} -1\right). 
\end{equation} 
For all positive integers $n$ and $k$, 
\begin{equation} \label{prim:recursionfk} 
\sum_{d=1}^n f_k\left( \left[\frac{n}{d}\right]\right) = {n\choose k} 
\end{equation} 
and 
\begin{equation} \label{prim:explicitfk} 
f_k(n) = \sum_{d=1}^{n} \mu(d) {[n/d] \choose k }. 
\end{equation} 
\et 

\noindent{\it Proof.} 
Let $A$ be a nonempty subset of $\{1,2,\ldots,n\}$. If $\gcd(A) = d$, then $A' = (1/d)\ast A = \{a/d : a\in A\}$ is a relatively prime subset of $\{1,2,\ldots,[n/d]\}$. Conversely, if $A'$ is a relatively prime subset of $\{1,2,\ldots,[n/d]\}$, 
then $A = d\ast A' = \{da' : a' \in A' \}$ is a nonempty subset of $\{1,2,\ldots,n\}$ with $\gcd(A)=d.$ It follows that there are exactly $f([n/d])$ subsets $A$ of $\{1,2,\ldots,n\}$ with $\gcd(A) = d$, and so 
\[ 
\sum_{d=1}^n f\left(\left[\frac{n}{d}\right]\right) = 2^n-1. 
\] 
We apply M\" obius inversion to the function $F(x) = f([x]).$ For all $x \geq 1$ we define 
\begin{align*} 
G(x) & = \sum_{1 \leq d \leq x} F\left(\frac{x}{d}\right) 
= \sum_{1 \leq d \leq x} f\left( \left[\frac{x}{d}\right] \right) 
= \sum_{d=1}^{[x]} f\left( \left[\frac{[x]}{d}\right] \right) 
= 2^{[x]} -1 
\end{align*} 
and so 
\begin{align*} 
f([x]) & = F(x) = \sum_{1 \leq d \leq x} \mu(d) G\left(\frac{x}{d}\right) = \sum_{d=1}^{[x]} \mu(d) \left(2^{[x/d]} -1\right). 
\end{align*} 
For $n \geq 1$ we have 
\[ 
f(n) = \sum_{d=1}^{n} \mu(d) \left(2^{[n/d]} -1\right). 
\] 

The proofs of~\eqref{prim:recursionfk} and~\eqref{prim:explicitfk} are similar. 
\hfill{$\Box$} 

\bt \label{prim:theorem:f(n)} 
For all positive integers $n$ and $k$, 
\[ 
2^n - 2^{[n/2]} - n2^{[n/3]} \leq f(n) \leq 2^n - 2^{[n/2]} 
\] 
and 
\[ 
{n \choose k} - {[n/2] \choose k} - n {[n/3] \choose k} 
\leq f_k(n) \leq {n \choose k} - {[n/2] \choose k}. 
\] 
\et 

\noindent{\it Proof.} 
For $n \geq 2$ we have 
\[ 
2^n = f(n) + f([n/2]) + \sum_{d=3}^n f\left(\left[\frac{n}{d}\right]\right) + 1 \leq f(n) + 2^{[n/2]} + n2^{[n/3]}. 
\] 
Combining this with~\eqref{prim:upper}, we obtain 
\[ 
2^n - 2^{[n/2]} - n2^{[n/3]} \leq f(n) \leq 2^n - 2^{[n/2]}. 
\] 
This also holds for $n=1$. 

The inequality for $f_k(n)$ follows similarly from~\eqref{prim:upper-fkn} and~\eqref{prim:recursionfk}. 
\hfill{$\Box$} 

Theorem~\ref{prim:theorem:f(n)} implies that $f(n) \sim 2^n$ as $n \rightarrow \infty,$ and so almost all finite sets of integers are relatively prime.

\section{A phi Function for Sets} 
The Euler phi function $\varphi(n)$ counts the number of positive integers $a \leq n$ such that $a$ is relatively prime to $n$. 
We define the function $\Phi(n)$ to be the number of nonempty subsets $A$ of $\{1,2,\ldots,n\}$ such that 
$\gcd(A)$ is relatively prime to $n$. 
For example, for distinct primes $p$ and $q$ we have 
\[ 
\Phi(p) = 2^p - 2 
\] 
\[ 
\Phi(p^2) = 2^{p^2} - 2^p 
\] 
and 
\[ 
\Phi(pq) = 2^{pq}- 2^q - 2^p +2. 
\] 
Define the function $\Phi_k(n)$ to be the number of subsets $A$ of $\{1,2,\ldots,n\}$ such that $\card(A)=k$ and 
$\gcd(A)$ is relatively prime to $n$. 
Note that $\Phi_1(n) =\varphi(n)$ for all $n \geq 1$.

\bt \label{prim:theorem:Phirecursion} 
For all positive integers $n$, 
\begin{equation} \label{prim:Phirecursion-1} 
\sum_{d|n} \Phi\left( d\right) = 2^n -1. 
\end{equation} 
Moreover, $\Phi(1)=1$ and, for $n \geq 2$, 
\begin{equation} \label{prim:Phirecursion-2} 
\Phi(n) = \sum_{d|n} \mu\left(d\right) 2^{n/d} 
\end{equation} 
where $\mu(n)$ is the M\" obius function. 
Similarly, for all positive integers $n$ and $k$, 
\begin{equation} \label{prim:Phirecursion-k1} 
\sum_{d|n} \Phi_k( d) = {n \choose k} 
\end{equation} 
and 
\begin{equation} \label{prim:Phirecursion-k2} 
\Phi_k(n) = \sum_{d|n} \mu\left(d\right) {n/d\choose k} 
\end{equation} 
\et 

\noindent{\it Proof.} 
For every divisor $d$ of $n$, we define the function $\Psi(n,d)$ to be the number of nonempty subsets $A$ of $\{1,2,\ldots,n\}$ such that the greatest common divisor of 
$\gcd(A)$ and $n$ is $d$. Thus, 
\[ 
\Psi(n,d) = 
\card\left( \left\{ A \subseteq \{1,2,\ldots,n\} :A \neq\emptyset \text{ and } \gcd(A \cup \{n\})=d \right\} \right). 
\] 
Then 
\[ 
\Psi(n,d) = \Phi\left( \frac{n}{d} \right) 
\] 
and 
\[ 
2^n -1 = \sum_{d|n} \Psi(n,d) = \sum_{d|n} \Phi\left( \frac{n}{d}\right) = \sum_{d|n} \Phi(d). 
\] 
We have $\Phi(1)=1.$ For $n \geq 2$ we apply the usual M\" obius inversion and obtain 
\begin{align*} 
\Phi(n) & = \sum_{d|n} \mu\left(d\right) \left( 2^{n/d} -1\right) \\ 
& = \sum_{d|n} \mu\left(d\right) 2^{n/d} 
- \sum_{d|n} \mu\left(d\right) \\ 
& = \sum_{d|n} \mu\left(d\right) 2^{n/d} 
\end{align*} 
since $\sum_{d|n}\mu(n/d)=0$ for $n \geq 2.$ 

The proofs of~\eqref{prim:Phirecursion-k1} and~\eqref{prim:Phirecursion-k2} are similar. 
\hfill{$\Box$}

\bt \label{prim:theorem:Phi} 
If $n$ is odd, then 
\[ 
\Phi(n) = 2^n + O\left( n2^{n/3} \right) 
\] 
and 
\[ 
\Phi_k(n) = {n \choose k} + O\left( n{[n/3] \choose k} \right). 
\] 
If $n$ is even, then 
\[ 
\Phi(n) = 2^n - 2^{n/2} + O\left( n2^{n/3} \right) 
\] 
and 
\[ 
\Phi_k(n) = {n \choose k} - {n/2 \choose k} + O\left( n{[n/3] \choose k} \right). 
\] 
\et 

\noindent{\it Proof.} 
We have 
\begin{align*} 
\Phi(n) 
& = \sum_{\substack{d=1\\ \gcd(d,n)=1}}^n \card\left( \left\{ A \subseteq \{1,2,\ldots,n\} :A \neq\emptyset \text{ and } \gcd(A)=d \right\} \right) \\ 
& = \sum_{\substack{d=1\\ \gcd(d,n)=1}}^n f([n/d]). 
\end{align*} 
Applying Theorem~\ref{prim:theorem:f(n)}, we see that if $n$ is odd, then 
\begin{align*} 
\Phi(n) 
& = f(n) + f([n/2]) + \sum_{ \substack{d=3\\ \gcd(d,n)=1} }^n f([n/d]) \\ 
& = \left( 2^n - 2^{[n/2]} + O\left(n2^{n/3}\right)\right) 
+ \left( 2^{[n/2]} + O\left(2^{n/4}\right)\right) 
+O\left( n2^{n/3} \right) \\ 
& = 2^n + O\left( n2^{n/3} \right). 
\end{align*} 
If $n$ is even, then 
\begin{align*} 
\Phi(n) 
& = f(n) + \sum_{ \substack{d=3\\ \gcd(d,n)=1} }^n f([n/d]) \\ 
& = \left( 2^n - 2^{n/2} + O\left(n2^{n/3}\right)\right) 
+O\left( n2^{n/3} \right) \\ 
& = 2^n - 2^{n/2} + O\left( n2^{n/3} \right). 
\end{align*} 
These estimates for $\Phi(n)$ also follow from identity~\eqref{prim:Phirecursion-2}. The estimates for $\Phi_k(n)$ follow from identity~\eqref{prim:Phirecursion-k2}. 
This completes the proof. 
\hfill{$\Box$} 

\vskip 30pt 
\noindent 
{\bf Acknowledgements.} I thank Greg Martin for the observation 
that~\eqref{prim:explicit} and~\eqref{prim:explicitfk} follow 
from~\eqref{prim:recursion} and~\eqref{prim:recursionfk} by M\" obius 
inversion, and Kevin O'Bryant for helpful discussions.

\providecommand{\bysame}{\leavevmode\hbox to3em{\hrulefill}\thinspace} 
\providecommand{\MR}{\relax\ifhmode\unskip\space\fi MR } 
\providecommand{\MRhref}[2]{%
\href{http://www.ams.org/mathscinet-getitem?mr=#1}{#2} 
} 
\providecommand{\href}[2]{#2}


\begin{thebibliography}{1} \footnotesize 

\bibitem{nath00A} 
M.~B. Nathanson, \emph{Elementary Methods in Number Theory}, Graduate Texts in Mathematics, vol. 195, Springer-Verlag, New York, 2000. 

\end{thebibliography}
\end{document}